\documentclass[10pt]{article}
\usepackage{cite}
\usepackage{mathrsfs}
\usepackage{amsfonts}
\usepackage{amsmath}
\usepackage{amsfonts,amssymb,color}
\usepackage{dsfont}
\usepackage{curves}
\usepackage{mathrsfs}
\usepackage{pifont}
\usepackage{amssymb}
\allowdisplaybreaks

\numberwithin{equation}{section}

\date{}

\def\BigRoman{\uppercase\expandafter{\romannumeral\number\count 255 }}
\def\Romannumeral{\afterassignment\BigRoman\count255=}

\begin{document}
\title{Distance spectral radius and $H_b$-factors in graphs
}
\author{\small  Jie Wu\footnote{Corresponding
author. E-mail address: jskjdxwj@126.com}\\
\small  School of Economics and management,\\
\small  Jiangsu University of Science and Technology,\\
\small  Zhenjiang, Jiangsu 212100, China\\
}

\maketitle
\begin{abstract}
\noindent Let $G$ be a connected graph, and let $b\geq2$ be an even integer. The distance spectral radius of $G$ is denoted by $\mu(G)$. An $H_b$-factor of $G$ is a spanning subgraph $F$ of $G$ with
$d_F(v)\in\{1,3,5,\ldots,b-1,b\}$ for any $v\in V(G)$, where $d_F(v)$ is the degree of $v$ in $F$. Lu and Wang provided a sufficient condition with respect to the number of odd components in $G-S$
for a connected graph $G$ of even order to contain an $H_b$-factor, where $S$ is a vertex subset of $G$ [H. Lu, D. Wang, On Cui-Kano's characterization problem on graph factors, J. Graph Theory 74
(2013) 335--343]. In this paper, motivated by Lu and Wang's above result, we establish an upper bound on the distance spectral radius $\mu(G)$ of a connected graph $G$ to guarantee that $G$ contains
an $H_b$-factor.
\\
\begin{flushleft}

(2020) Mathematics Subject Classification: 05C50, 05C70, 90B99
\end{flushleft}
\end{abstract}

\section{Introduction}

Graphs considered in this paper are finite, undirected and simple. Let $G$ be a graph with vertex set $V(G)$ and edge set $E(G)$. We denote the order of $G$ by $|V(G)|=n$. For a vertex $v\in V(G)$,
the degree of $v$ in $G$ is denoted by $d_G(v)$. Let $\delta(G)$ and $o(G)$ stand for the minimum degree and the number of odd components in $G$, respectively. For a subset $S\subseteq V(G)$, the
subgraphs of $G$ induced by $S$ and $V(G)\setminus S$ are denoted by $G[S]$ and $G-S$, respectively. Let $K_n$ denote the complete graph with $n$ vertices. Given two vertex-disjoint graphs $G_1$ and
$G_2$, the union and the join of $G_1$ and $G_2$ are denoted by $G_1\cup G_2$ and $G_1\vee G_2$, respectively.

For a given connected graph $G$ with $V(G)=\{v_1,v_2,\ldots,v_n\}$, the adjacency matrix of $G$, denoted by $\mathcal{D}(G)$, is defined as $\mathcal{D}(G)=(d_{ij})_{n\times n}$, where $d_{ij}$ equals
the length of a shortest path from $v_i$ to $v_j$ in $G$, which is called the distance between $v_i$ and $v_j$ in $G$. The largest eigenvalue of $\mathcal{D}(G)$, denoted by $\mu(G)$, is called the
distance spectral radius of $G$. Some properties of the distance spectral radius in $G$ can be found in \cite{NP,Em,MK,ZLi,Zs,ZBS,Zt}.

Let $b$ be a positive even integer. Then an odd $[1,b-1]$-factor of $G$ is a spanning subgraph $F$ of $G$ with $d_F(v)\in\{1,3,5,\ldots,b-1\}$ for any $v\in V(G)$. In particular, an odd $[1,1]$-factor
is called a perfect matching. An $H_b$-factor of $G$ is a spanning subgraph $F$ of $G$ with $d_F(v)\in\{1,3,5,\ldots,b-1,b\}$ for any $v\in V(G)$.

Zhang and Lin \cite{ZL} established some relationships between the distance spectral radius and perfect matchings in graphs and bipartite graphs. Li and Miao \cite{LM} arose a distance spectral radius
condition for the existence of an odd $[1,b-1]$-factor in a graph. Hu and Zhang \cite{HZ} proposed a sufficient condition based on the distance spectral radius which guarantees the existence of an odd
$[1,b-1]$-factor in a 1-binding graph. Fan, Lin and Lu \cite{FLL}, Zhou \cite{Zs1} provided some sufficient conditions for a graph to possess an odd $[1,b-1]$-factor. Lu and Wang \cite{LW} presented
some sufficient conditions for the existence of an $H_b$-factor in a connected graph. Fan and Liu \cite{FL} showed a sufficient condition in terms of the adjacency spectral radius to guarantee that a
connected graph contains an $H_b$-factor. Zhou \cite{Zs2} gave two sufficient conditions based on the size and the signless Laplacian spectral radius to ensure the existence of an $H_b$-factor in a
connected graph. Some other results on graph factors were obtained by Wu \cite{Wc,Ws,Wu}, Hartvigsen \cite{Hf}, Chiba, Egawa and Fujita \cite{CEF}, Zhou et al \cite{Zs3,ZBW}, Pan and Zhou \cite{PZ},
Wu, Zhou and Liu \cite{WZL}.

Motivated by \cite{LW} directly, we pose a distance spectral radius condition for the existence of an $H_b$-factor in a graph. Our main result is shown in the following.

\medskip

\noindent{\textbf{Theorem 1.1.}} Let $G$ be a 2-connected graph of even order $n\geq\max\{(b+5)\delta+b,\frac{2}{3}b\delta^{3}+\delta^{2}-3b\delta\}$ with minimum degree $\delta\geq2$, where $b\geq2$ is
an even integer. If
$$
\mu(G)\leq\mu(K_{\delta}\vee(K_{n-(b+1)\delta}\cup b\delta K_1)),
$$
then $G$ contains an $H_b$-factor unless $G=K_{\delta}\vee(K_{n-(b+1)\delta}\cup b\delta K_1)$.

\medskip

\section{Tools}

Lu and Wang \cite{LW} presented a sufficient condition for a connected graph with even order to possess an $H_b$-factor.

\medskip

\noindent{\textbf{Lemma 2.1}} (Lu and Wang \cite{LW}). Let $b$ be a positive even integer, and let $G$ be a connected graph of even order. If
$$
o(G-S)\leq b|S|
$$
for any nonempty subset $S\subseteq V(G)$, then $G$ contains an $H_b$-factor.

\medskip

\noindent{\textbf{Lemma 2.2}} (Minc \cite{Mn}). Let $G$ be a connected graph with two nonadjacent vertices $u,v\in V(G)$. Then
$$
\mu(G)>\mu(G+uv),
$$
where $G+uv$ denotes the graph obtained from $G$ by adding an edge to connect $u$ and $v$.

\medskip

We denote the Wiener index of a connected graph $G$ with $n$ vertices by $W(G)=\sum\limits_{i<j}d_{ij}$. The following lemma is easily obtained based on the Rayleigh quotient \cite{HJ}.

\medskip

\noindent{\textbf{Lemma 2.3.}} Let $G$ be a connected graph of order $n$. Then
$$
\mu(G)=\max_{X\neq\mathbf{0}}\frac{X^{T}\mathcal{D}(G)X}{X^{T}X}\geq\frac{\mathbf{1}^{T}\mathcal{D}(G)\mathbf{1}}{\mathbf{1}^{T}\mathbf{1}}=\frac{2W(G)}{n},
$$
where $\mathbf{1}=(1,1,\ldots,1)^{T}$.

\medskip

\noindent{\textbf{Lemma 2.4}} (Zheng, Li, Luo and Wang \cite{ZLLW}). Let $s+\sum\limits_{i=1}^{t}n_i=n$. If $n_1\geq n_2\geq\cdots\geq n_t\geq p\geq1$ and $n_1\leq n-s-p(t-1)$, then
$$
\mu(K_s\vee(K_{n_1}\cup K_{n_2}\cup\cdots\cup K_{n_t}))\geq\mu(K_s\vee(K_{n-s-p(t-1)}\cup(t-1)K_p)),
$$
with equality if and only if $(n_1,n_2,\ldots,n_t)=(n-s-p(t-1),p,\ldots,p)$.

\medskip

Assume that $M$ be a real symmetric matrix whose rows and columns are indexed by $V=\{1,2,\ldots,n\}$. Suppose that $M$, based on the partition $\pi: V=V_1\cup V_2\cup\cdots\cup V_r$, can be written as
\begin{align*}
M=\left(
  \begin{array}{ccc}
    M_{11} & \cdots & M_{1,r}\\
    \vdots & \ddots & \vdots\\
    M_{r1} & \cdots & M_{r,r}\\
  \end{array}
\right),
\end{align*}
where $M_{ij}$ is the submatrix of $M$ derived by rows in $V_i$ and columns in $V_j$. Let $q_{ij}$ equal the average row sum of $M_{ij}$. Then the matrix $M_{\pi}=(q_{ij})_{r\times r}$ is said to be the
quotient matrix of $M$ with respect to the partition $\pi$. If the row sum of every block $M_{ij}$ is a constant for $1\leq i,j\leq r$, then the partition $\pi$ is equitable.

\medskip

\noindent{\textbf{Lemma 2.5}} (You, Yang, So and Xi \cite{YYSX}). Let $M$ be a real nonnegative matrix with an equitable partition $\pi$, and let $M_{\pi}$ be the quotient matrix of $M$ based on the partition
$\pi$. Then the largest eigenvalues of $M$ and $M_{\pi}$ are equal.

\medskip

\section{The proof of Theorem 1.1}

\noindent{\it Proof of Theorem 1.1.} Suppose that $G$ contains no $H_b$-factor. Then from Lemma 2.1, there exists a nonempty subset $S$ of $V(G)$ such that $o(G-S)\geq b|S|+1$. Let $|S|=s$. Then $G$ is a
spanning subgraph of $G_1=K_s\vee(K_{n_1}\cup K_{n_2}\cup\cdots\cup K_{n_{bs+1}})$ for some positive odd integers $n_1\geq n_2\geq\cdots\geq n_{bs+1}$ with $\sum\limits_{i=1}^{bs+1}n_i=n-s$. According to
Lemma 2.2, we have
\begin{align}\label{eq:3.1}
\mu(G)\geq\mu(G_1),
\end{align}
where the equality occurs if and only if $G=G_1$. In what follows, we divide the proof into three cases.

\noindent{\bf Case 1.} $s\geq\delta+1$.

Let $G_2=K_s\vee(K_{n-(b+1)s}\cup bsK_1)$, where $n\geq(b+1)s+1$. In view of Lemma 2.4, we obtain
\begin{align}\label{eq:3.2}
\mu(G_1)\geq\mu(G_2),
\end{align}
with equality holding if and only if $(n_1,n_2,\ldots,n_{bs+1})=(n-(b+1)s,1,\ldots,1)$.

Based on the partition $V(G_2)=V(K_s)\cup V(K_{n-(b+1)s})\cup V(bsK_1)$, the quotient matrix of $\mathcal{D}(G_2)$ is equal to
\begin{align*}
B_2=\left(
  \begin{array}{ccc}
    s-1 & n-(b+1)s & bs\\
    s & n-(b+1)s-1 & 2bs\\
    s & 2n-2(b+1)s & 2bs-2\\
  \end{array}
\right).
\end{align*}
Then the characteristic polynomial of $B_2$ is
\begin{align*}
\phi_{B_2}(x)=&x^{3}-(n+bs-4)x^{2}-(2bsn+3n-(2b^{2}+3b)s^{2}+bs-5)x\\
&+(bs^{2}-2bs-2)n-(b^{2}+b)s^{3}+(2b^{2}+3b)s^{2}+2.
\end{align*}
In terms of Lemma 2.5 and the equitable partition $V(G_2)=V(K_s)\cup V(K_{n-(b+1)s})\cup V(bsK_1)$, $\mu(G_2)$ is equal to the largest root of $\phi_{B_2}(x)=0$.

Let $G_*=K_{\delta}\vee(K_{n-(b+1)\delta}\cup b\delta K_1)$. Given the partition $V(G_*)=V(K_{\delta})\cup V(K_{n-(b+1)\delta})\cup V(b\delta K_1)$, the quotient matrix of $\mathcal{D}(G_*)$ is
\begin{align*}
B_*=\left(
  \begin{array}{ccc}
    \delta-1 & n-(b+1)\delta & b\delta\\
    \delta & n-(b+1)\delta-1 & 2b\delta\\
    \delta & 2n-2(b+1)\delta & 2b\delta-2\\
  \end{array}
\right),
\end{align*}
and the characteristic polynomial of $B_*$ is given by
\begin{align*}
\phi_{B_*}(x)=&x^{3}-(n+b\delta-4)x^{2}-(2b\delta n+3n-(2b^{2}+3b)\delta^{2}+b\delta-5)x\\
&+(b\delta^{2}-2b\delta-2)n-(b^{2}+b)\delta^{3}+(2b^{2}+3b)\delta^{2}+2.
\end{align*}
Note that the partition $V(G_*)=V(K_{\delta})\cup V(K_{n-(b+1)\delta})\cup V(b\delta K_1)$ is equitable. According to Lemma 2.5, $\mu(G_*)$ equals the largest root of $\phi_{B_*}(x)=0$.

Write $\mu(G_*)=\mu$. By virtue of $n\geq(b+5)\delta+b$, we get
\begin{align}\label{eq:3.3}
\mu\geq&\frac{2W(G_*)}{n}\nonumber\\
=&\frac{n^{2}+(2b\delta-1)n-(b^{2}+2b)\delta^{2}-b\delta}{n}\nonumber\\
>&n+b\delta-1.
\end{align}

Notice that $\phi_{B_*}(\mu)=0$. Then
\begin{align}\label{eq:3.4}
\phi_{B_2}(\mu)=&\phi_{B_2}(\mu)-\phi_{B_*}(\mu)\nonumber\\
=&(s-\delta)\Big(-b\mu^{2}+(-2bn+(2b^{2}+3b)s+(2b^{2}+3b)\delta-b)\mu\nonumber\\
&+(bs+b\delta-2b)n-(b^{2}+b)s^{2}+(2b^{2}+3b-b^{2}\delta-b\delta)s\nonumber\\
&-(b^{2}+b)\delta^{2}+(2b^{2}+3b)\delta\Big).
\end{align}
Let $f(x)=-bx^{2}+(-2bn+(2b^{2}+3b)s+(2b^{2}+3b)\delta-b)x+(bs+b\delta-2b)n-(b^{2}+b)s^{2}+(2b^{2}+3b-b^{2}\delta-b\delta)s-(b^{2}+b)\delta^{2}+(2b^{2}+3b)\delta$. By a simple calculation, the derivative
function of $f(x)$ is
$$
f'(x)=-2bx-2bn+(2b^{2}+3b)s+(2b^{2}+3b)\delta-b.
$$
When $x\geq n+b\delta-1$, it follows from $s\geq\delta+1$ and $n\geq(b+1)s+1$ that
\begin{align*}
f'(x)\leq&-2b(n+b\delta-1)-2bn+(2b^{2}+3b)s+(2b^{2}+3b)\delta-b\\
=&-4bn+(2b^{2}+3b)s+3b\delta+b\\
\leq&-4b((b+1)s+1)+(2b^{2}+3b)s+3b\delta+b\\
=&-(2b^{2}+b)s+3b\delta-3b\\
\leq&-(2b^{2}+b)(\delta+1)+3b\delta-3b\\
=&-(2b^{2}-2b)\delta-2b^{2}-4b\\
<&0,
\end{align*}
which implies that $f(x)$ is decreasing for $x\geq n+b\delta-1$. Together with \eqref{eq:3.3}, we obtain
\begin{align*}
f(\mu)\leq&f(n+b\delta-1)\\
=&-b(n+b\delta-1)^{2}+(-2bn+(2b^{2}+3b)s+(2b^{2}+3b)\delta-b)(n+b\delta-1)\\
&+(bs+b\delta-2b)n-(b^{2}+b)s^{2}+(2b^{2}+3b-b^{2}\delta-b\delta)s\\
&-(b^{2}+b)\delta^{2}+(2b^{2}+3b)\delta\\
=&-(b^{2}+b)s^{2}+((2b^{2}+4b)n+(2b^{3}+2b^{2}-b)\delta)s\\
&-3bn^{2}+((-2b^{2}+4b)\delta+b)n+(b^{3}+2b^{2}-b)\delta^{2}+b^{2}\delta\\
\leq&-(b^{2}+b)\left(\frac{n-1}{b+1}\right)^{2}+((2b^{2}+4b)n+(2b^{3}+2b^{2}-b)\delta)\left(\frac{n-1}{b+1}\right)\\
&-3bn^{2}+((-2b^{2}+4b)\delta+b)n+(b^{3}+2b^{2}-b)\delta^{2}+b^{2}\delta \ \ \ \ \ \left(\mbox{since} \ s\leq\frac{n-1}{b+1}\right)\\
=&\frac{1}{b+1}\Big(-b^{2}n^{2}+((4b^{2}+3b)\delta-b^{2}-b)n+(b^{4}+3b^{3}+b^{2}-b)\delta^{2}-(b^{3}+b^{2}-b)\delta-b\Big)\\
\leq&\frac{1}{b+1}\Big(-b^{2}((b+5)\delta+b)^{2}+((4b^{2}+3b)\delta-b^{2}-b)((b+5)\delta+b)\\
&+(b^{4}+3b^{3}+b^{2}-b)\delta^{2}-(b^{3}+b^{2}-b)\delta-b\Big) \ \ \ \ \ (\mbox{since} \ n\geq(b+5)\delta+b)\\
=&\frac{1}{b+1}\Big((-3b^{3}-b^{2}+14b)\delta^{2}+(-2b^{4}-8b^{3}-4b^{2}+b)\delta-b^{4}-b^{3}-b^{2}-b\Big)\\
<&0 \ \ \ \ \ (\mbox{since} \ \delta\geq2 \ \mbox{and} \ b\geq2).
\end{align*}
Combining this with \eqref{eq:3.4} and $s\geq\delta+1$, we conclude $\phi_{B_2}(\mu)=(s-\delta)f(\mu)<0$. This yields $\mu(G_*)=\mu<\mu(G_2)$. Together with \eqref{eq:3.1} and \eqref{eq:3.2}, we deduce
$$
\mu(G)\geq\mu(G_1)\geq\mu(G_2)>\mu(G_*)=\mu(K_{\delta}\vee(K_{n-(b+1)\delta}\cup b\delta K_1)),
$$
which contradicts $\mu(G)\leq\mu(K_{\delta}\vee(K_{n-(b+1)\delta}\cup b\delta K_1))$.

\noindent{\bf Case 2.} $s=\delta$.

Recall that $G_1=K_s\vee(K_{n_1}\cup K_{n_2}\cup\cdots\cup K_{n_{bs+1}})$. Based on $s=\delta$ and Lemma 2.4, we get
$$
\mu(G_1)\geq\mu(K_{\delta}\vee(K_{n-(b+1)\delta}\cup b\delta K_1)),
$$
with equality holding if and only if $G_1=K_{\delta}\vee(K_{n-(b+1)\delta}\cup b\delta K_1)$. Together with \eqref{eq:3.1}, we have
$$
\mu(G)\geq\mu(K_{\delta}\vee(K_{n-(b+1)\delta}\cup b\delta K_1)),
$$
with equality occurring if and only if $G=K_{\delta}\vee(K_{n-(b+1)\delta}\cup b\delta K_1)$, a contradiction.

\noindent{\bf Case 3.} $s\leq\delta-1$.

Let $G_3=K_s\vee(K_{n-s-bs(\delta+1-s)}\cup bsK_{\delta+1-s})$. Recall that $G$ is a spanning subgraph of $G_1=K_s\vee(K_{n_1}\cup K_{n_2}\cup\cdots\cup K_{n_{bs+1}})$ for some positive odd integers
$n_1\geq n_2\geq\cdots\geq n_{bs+1}$ with $\sum\limits_{i=1}^{bs+1}n_i=n-s$. Since $\delta(G_1)\geq\delta(G)=\delta$, we deduce $n_{bs+1}\geq\delta+1-s$. According to Lemma 2.4, we conclude
\begin{align}\label{eq:3.5}
\mu(G_1)\geq\mu(G_3),
\end{align}
where the equality follows if and only if $(n_1,n_2,\ldots,n_{bs+1})=(n-s-bs(\delta+1-s),\delta+1-s,\ldots,\delta+1-s)$. Notice that $n$ is even and $G$ is 2-connected. If $s=1$, then
$1=o(G-S)\geq bs+1=b+1>1$, which is impossible. Next, we shall consider $2\leq s\leq\delta-1$.

In terms of the partition $V(G_3)=V(K_s)\cup V(K_{n-s-bs(\delta+1-s)})\cup V(bsK_{\delta+1-s})$, the quotient matrix of $\mathcal{D}(G_3)$ is given by
\begin{align*}
B_3=\left(
  \begin{array}{ccc}
    s-1 & n-s-bs(\delta+1-s) & b(\delta+1-s)\\
    s & n-s-bs(\delta+1-s)-1 & 2bs(\delta+1-s)\\
    s & 2n-2s-2bs(\delta+1-s) & 2(bs-1)(\delta+1-s)+\delta-s\\
  \end{array}
\right),
\end{align*}
and the characteristic polynomial of $B_3$ equals
\begin{align*}
\phi_{B_3}(x)=&x^{3}+(-n-bs(\delta+1-s)+\delta-s+4)x^{2}\\
&+\Big(2bs^{4}-(4b^{2}\delta+4b^{2}+2b)s^{3}+(2bn+2b^{2}\delta^{2}+4b^{2}\delta+b\delta+2b^{2}+3b)s^{2}\\
&-(2b\delta n+2bn-n-b\delta^{2}+b+2)s-\delta n-3n+2\delta+5\Big)x\\
&-b^{2}s^{5}+(2b^{2}\delta+4b^{2}+b)s^{4}-(bn+b^{2}\delta^{2}+6b^{2}\delta+b\delta+5b^{2}+3b)s^{3}\\
&+(b\delta n+3bn+2b^{2}\delta^{2}+4b^{2}\delta+b\delta+2b^{2}+2b)s^{2}\\
&-(2b\delta n+2bn-n-b\delta^{2}-b\delta+1)s-\delta n-2n+\delta+2.
\end{align*}
By means of Lemma 2.5 and the equitable partition $V(G_3)=V(K_s)\cup V(K_{n-s-bs(\delta+1-s)})\cup V(bsK_{\delta+1-s})$, the largest root of $\phi_{B_3}(x)=0$ equals $\mu(G_3)$.

Recall that $\phi_{B_*}(\mu)=0$. Then
\begin{align}\label{eq:3.6}
\phi_{B_3}(\mu)=&\phi_{B_3}(\mu)-\phi_{B_*}(\mu)=(\delta-s)g(\mu),
\end{align}
where $g(x)=(-bs+b+1)x^{2}+(-2b^{2}s^{3}+(2b^{2}\delta+4b^{2}+2b)s^{2}+(-2bn+b\delta-2b^{2}-3b)s+2bn-n-2b^{2}\delta-3b\delta+b+2)x+b^{2}s^{4}-(b^{2}\delta+4b^{2}+b)s^{3}+(bn+2b^{2}\delta+5b^{2}+3b)s^{2}
+(-3bn+b^{2}\delta+2b\delta-2b^{2}-2b)s-b\delta n+2bn-n+b^{2}\delta^{2}+b\delta^{2}-2b^{2}\delta-3b\delta+1$. For $x\geq n+b\delta-1$, it follows from $b\geq2$, $2\leq s\leq\delta-1$ and
$n\geq\frac{2}{3}b\delta^{3}+\delta^{2}-3b\delta>\frac{2b^{2}\delta^{3}+2b\delta^{2}-(10b^{2}+3b)\delta-16}{4b-1}$ that
\begin{align*}
g'(x)=&2(-bs+b+1)x-2b^{2}s^{3}+(2b^{2}\delta+4b^{2}+2b)s^{2}\\
&+(-2bn+b\delta-2b^{2}-3b)s+2bn-n-2b^{2}\delta-3b\delta+b+2\\
\leq&2(-bs+b+1)(n+b\delta-1)-2b^{2}s^{3}+(2b^{2}\delta+4b^{2}+2b)s^{2}\\
&+(-2bn+b\delta-2b^{2}-3b)s+2bn-n-2b^{2}\delta-3b\delta+b+2\\
=&-2b^{2}s^{3}+(2b^{2}\delta+4b^{2}+2b)s^{2}-(4bn+2b^{2}\delta-b\delta+2b^{2}+b)s+(4b+1)n-b\delta-b\\
\leq&-16b^{2}+(2b^{2}\delta+4b^{2}+2b)(\delta-1)^{2}-2(4bn+2b^{2}\delta-b\delta+2b^{2}+b)+(4b+1)n-b\delta-b\\
=&-(4b-1)n+2b^{2}\delta^{3}+2b\delta^{2}-(10b^{2}+3b)\delta-16\\
<&0,
\end{align*}
which implies that $g(x)$ is decreasing in the interval $[n+b\delta-1,+\infty)$. Combining this with \eqref{eq:3.3}, $b\geq2$ and $2\leq s\leq\delta-1$, we deduce
\begin{align}\label{eq:3.7}
g(\mu)<&g(n+b\delta-1)\nonumber\\
=&(-bs+b+1)(n+b\delta-1)^{2}+\Big(-2b^{2}s^{3}+(2b^{2}\delta+4b^{2}+2b)s^{2}\nonumber\\
&+(-2bn+b\delta-2b^{2}-3b)s+2bn-n-2b^{2}\delta-3b\delta+b+2\Big)(n+b\delta-1)+b^{2}s^{4}\nonumber\\
&-(b^{2}\delta+4b^{2}+b)s^{3}+(bn+2b^{2}\delta+5b^{2}+3b)s^{2}+(-3bn+b^{2}\delta+2b\delta-2b^{2}-2b)s\nonumber\\
&-b\delta n+2bn-n+b^{2}\delta^{2}+b\delta^{2}-2b^{2}\delta-3b\delta+1\nonumber\\
=&b^{2}s^{4}-(2b^{2}n+2b^{3}\delta+b^{2}\delta+2b^{2}+b)s^{3}\nonumber\\
&+((2b^{2}\delta+4b^{2}+3b)n+2b^{3}\delta^{2}+4b^{3}\delta+2b^{2}\delta+b^{2}+b)s^{2}\nonumber\\
&-(3bn^{2}+(4b^{2}\delta-b\delta+2b^{2}+2b)n+b^{3}\delta^{2}-b^{2}\delta^{2}+2b^{3}\delta-b\delta)s\nonumber\\
&+3bn^{2}+(2b^{2}\delta-3b\delta-b)n-b^{3}\delta^{2}-b^{2}\delta^{2}+b\delta^{2}-b^{2}\delta\nonumber\\
\leq&b^{2}(\delta-1)^{4}-8(2b^{2}n+2b^{3}\delta+b^{2}\delta+2b^{2}+b)\nonumber\\
&+((2b^{2}\delta+4b^{2}+3b)n+2b^{3}\delta^{2}+4b^{3}\delta+2b^{2}\delta+b^{2}+b)(\delta-1)^{2}\nonumber\\
&-2(3bn^{2}+(4b^{2}\delta-b\delta+2b^{2}+2b)n+b^{3}\delta^{2}-b^{2}\delta^{2}+2b^{3}\delta-b\delta)\nonumber\\
&+3bn^{2}+(2b^{2}\delta-3b\delta-b)n-b^{3}\delta^{2}-b^{2}\delta^{2}+b\delta^{2}-b^{2}\delta\nonumber\\
=&-3bn^{2}+(2b^{2}\delta^{3}+3b\delta^{2}-12b^{2}\delta-7b\delta-16b^{2}-2b)n\nonumber\\
&+2b^{3}\delta^{4}+b^{2}\delta^{4}-2b^{2}\delta^{3}-9b^{3}\delta^{2}+4b^{2}\delta^{2}+2b\delta^{2}-16b^{3}\delta-13b^{2}\delta-14b^{2}-7b\nonumber\\
\triangleq&h(n).
\end{align}
Note that
$$
\frac{2b^{2}\delta^{3}+3b\delta^{2}-12b^{2}\delta-7b\delta-16b^{2}-2b}{6b}<\frac{2}{3}b\delta^{3}+\delta^{2}-3b\delta\leq n
$$
by $\delta\geq s+1\geq3$ and the condition of the theorem. Then
\begin{align*}
h(n)\leq&h\Big(\frac{2}{3}b\delta^{3}+\delta^{2}-3b\delta\Big)\\
=&-3b\Big(\frac{2}{3}b\delta^{3}+\delta^{2}-3b\delta\Big)^{2}+(2b^{2}\delta^{3}+3b\delta^{2}-12b^{2}\delta-7b\delta-16b^{2}-2b)\Big(\frac{2}{3}b\delta^{3}+\delta^{2}-3b\delta\Big)\\
&+2b^{3}\delta^{4}+b^{2}\delta^{4}-2b^{2}\delta^{3}-9b^{3}\delta^{2}+4b^{2}\delta^{2}+2b\delta^{2}-16b^{3}\delta-13b^{2}\delta-14b^{2}-7b\\
=&-\frac{11}{3}b^{2}\delta^{4}-\frac{32}{3}b^{3}\delta^{3}-\frac{19}{3}b^{2}\delta^{3}-7b\delta^{3}+9b^{2}\delta^{2}+32b^{3}\delta-7b^{2}\delta-14b^{2}-7b\\
<&0 \ \ \ \ \ (\mbox{since} \ b\geq2 \ \mbox{and} \ \delta\geq3).
\end{align*}
Combining this with \eqref{eq:3.6}, \eqref{eq:3.7} and $2\leq s\leq\delta-1$, we conclude
$$
\phi_{B_3}(\mu)=(\delta-s)g(\mu)<(\delta-s)h(n)<0.
$$
This implies that $\mu(G_3)>\mu=\mu(G_*)$. Together with \eqref{eq:3.1} and \eqref{eq:3.5}, we infer
$$
\mu(G)\geq\mu(G_1)\geq\mu(G_3)>\mu(G_*)=\mu(K_{\delta}\vee(K_{n-(b+1)\delta}\cup b\delta K_1)),
$$
which contradicts $\mu(G)\leq\mu(K_{\delta}\vee(K_{n-(b+1)\delta}\cup b\delta K_1))$. Theorem 1.1 is proved. \hfill $\Box$

\medskip

\section*{Declaration of competing interest}

The author declares that he has no known competing financial interests or personal relationships that could have appeared to influence the work reported in this paper.

\section*{Data availability}

No data was used for the research described in the article.

\medskip



\end{document}